\documentclass{article}
\usepackage{amssymb,amsmath,enumerate,latexsym,geometry}
\usepackage[all]{xy}
\begin{document}

\begin{center}
\textbf{\LARGE Rewriting as a Special Case of }\\
\textbf{\LARGE Noncommutative Gr\"obner Basis Theory}\\
{\Large Anne Heyworth}\\
{\Large University of Wales, Bangor}
\end{center}
 
\section{Introduction}

Rewriting for semigroups is a special case of Gr\"obner basis theory for
noncommutative polynomial algebras. The fact is a kind of folklore but
is not fully recognised. So our aim in this paper is to elucidate this
relationship.\\

A good introduction to string rewriting is \cite{BoOt}, and a recent
introduction to noncommutative Gr\"obner basis theory is \cite{Ufn}.
Similarities between the two critical pair completion methods (Knuth-
Bendix  and Buchberger's algorithm) have often been pointed out in the
commutative case. The connection was first observed in
\cite{Loos,BuLoos} and more closely analysed in \cite{Bu86,Bu87} and
more recently in \cite{Sto} and \cite{Birgitsthesis}.
In particular it is well known that the commutative Buchberger algorithm
may be applied to presentations of abelian groups to obtain complete
rewrite systems.\\

Rewriting involves a presentation $sgp \langle X | R \rangle$ of a semigroup
$S$ and presents $S$ as a factor semigroup $X^\dagger / =_R$ where
$X^\dagger$ is the free semigroup on $X$ and $=_R$ is the congruence generated
by the subset $R$ of $X^\dagger \times X^\dagger$. Noncommutative
Gr\"obner basis theory involves a presentation $alg \langle X | F
\rangle$ of a noncommutative algebra $A$ over a field $K$ and
presents $A$ as a factor algebra $K[X^\dagger]/ \langle F \rangle$ where
$K[X^\dagger]$ is the free $K$-algebra on the semigroup $X^\dagger$
 and $\langle F \rangle$ is the ideal
generated by $F$, a subset of $K[X^\dagger]$.
Given a semigroup presentation $sgp\langle X | R \rangle$ we consider the
algebra presentation $alg \langle X | F \rangle$ where $F:=\{l-r : (l,r) \in R
\}$. 
It is well known that the word problem for $sgp \langle X | R \rangle$ is
solvable if and only if the (monomial) equality problem for $alg \langle X
| F \rangle$ is solvable. 
Teo Mora \cite{TMora} recorded that a complete rewrite system 
for a semigroup $S$ presented by $sgp\langle X|Rel\rangle$ is
equivalent to a noncommutative Gr\"obner basis for the ideal specified
by the congruence $=_R$ on $X^\dagger$ in the algebra
$\mathbb{F}_3[X^\dagger]$ where $\mathbb{F}_3$ is the field with
elements $\{-1,0,1\}$.\\ 

In this paper we show that the noncommutative Buchberger
algorithm applied to $F$ corresponds step-by-step to the Knuth-Bendix
completion procedure for $R$. This is the meaning intended for the first
sentence of this paper.

\section{Results}

First we note that the relation between the two kinds of presentation is
given by the following variation of a result of \cite{TMora}.\\

\textbf{\large Proposition}\\
Let $K$ be a field and let $S$ be a semigroup with presentation $sgp
\langle X | R \rangle$.
Then the algebra $K[S]$ is isomorphic to the factor algebra
$K[X^\dagger]/\! \langle F \rangle$ where $F$ is the basis $\{l-r|
(l,r) \in R\}$.

\textbf{Proof}\\
Define  $\phi:K[X^\dagger] \to K[S]$ by $\phi(k_1w_1 + \cdots + k_tw_t)
:= k_1[w_1]_R+\cdots+k_t[w_t]_R$ for $k_1, \ldots, k_t \in K$, $w_1,
\ldots, w_t
\in X^\dagger$.
Define a homomorphism $\phi':K[X^\dagger]/\! \langle F \rangle \to K[S]$ by
$\phi'([p]_F):=\phi(p)$.
It is injective since $\phi'[p]_F = \phi[q]_F$ if and only if
$[p]_F = [q]_F$ (using
the definitions $\phi(p)=\phi(q) \Leftrightarrow p-q \in \langle F \rangle$).
It is also surjective. Let $f\in K[S]$. Then $f = k_1m_1 + \cdots +
k_tm_t$
for some
$k_1,\ldots,k_t\in K$, $m_1,\ldots,m_t\in S$. Since $S$ is presented by
$sgp\langle X | R \rangle$ there exist $w_1, \ldots, w_t \in X^\dagger$ such
that
$[w_i]_R = m_i$ for
$i=1,\ldots,t$. Therefore let $p = k_1w_1 + \cdots + k_tw_t$. Clearly $p
\in K[X^\dagger]$ and also $\phi'[p]_F = f$. Hence $\phi'$ is an
isomorphism.{\hfill $\Box$\\}

Now we give our main result.\\

\textbf{\large Theorem}\\
Let $sgp\langle X | R\rangle$ be a semigroup presentation, let $K$ be a field
and let $alg\langle X | F \rangle$ be the $K$-algebra presentation with
$F:=\{ l-r: (l,r) \in R \}$.
Then the Knuth-Bendix completion algorithm for the rewrite system $R$ 
corresponds step-by-step to the noncommutative Buchberger algorithm for
finding a Gr\"obner basis for the ideal generated by $F$.\\

\textbf{Proof}
Both the Knuth-Bendix algorithm for $R$ and the Buchberger algorithm for
$F$ begin by specifying a monomial ordering on $X^\dagger$ which we
denote $>$.\\

The correspondence between terminology in the two cases is
\begin{alignat*}{2}
(i)    &\quad \text{rewrite system} &&\quad \text{basis} \\
(ii)   &\quad \text{rule }  &&\quad \text{two-term polynomial }  \\
(iii)  &\quad \text{word }  &&\quad \text{monomial }  \\
(iv)   &\quad \text{reduction}  &&\quad \text{reduction }  \\
(v)   &\quad \text{left hand side} &&\quad \text{leading monomial} \\
(vi)    &\quad \text{subword} &&\quad \text{submonomial} \\
(vii)   &\quad \text{right hand side} &&\quad \text{remainder} \\
(viii)  &\quad \text{overlap}  &&\quad \text{match}\\
(ix) &\quad \text{critical pair }  &&\quad \text{S-polynomial }
\intertext{This key part of the correspondence (viii) and (ix) is illustrated 
diagrammatically in the next section}
(x)    &\quad \text{resolve} &&\quad \text{reduce to zero} \\
(xi)   &\quad \text{reduced critical pair} &&\quad \text{reduced S-polynomial} \\
(xii)  &\quad \text{complete rewrite system} &&\quad \text{Gr\"obner basis}
\end{alignat*}

In terms of rewriting we consider the {rewrite system} $R$
which consists of a {set of rules} of the form $(l,r)$
orientated so that $l>r$.
A {word} $w \in X^\dagger$ may be {reduced} 
with
respect to $R$ if it contains the {left hand side} $l$ of a
rule $(l,r)$ as a {subword} i.e. if $w=ulv$ for some $u,v
\in
X^*$.
To reduce $w=ulv$ using the rule $(l,r)$ we replace $l$ by the
{right hand side} $r$ of the rule, and write $ulv \to_R
urv$.
The Knuth-Bendix algorithm looks for {overlaps between rules}.
Given a pair of rules $(l_1,r_1)$, $(l_2,r_2)$ there are four possible
ways
in which an overlap can occur: $l_1=u_2l_2v_2$, $u_1l_1v_1=l_2$,
$l_1v_1=u_2l_2$ and $u_1l_1=l_2v_2$. The {critical pair}
resulting from an overlap is the pair of words resulting from applying
each rule to the smallest word on which the overlap occurs. The critical
pairs resulting from each of the four overlaps are: $(r_1,u_2r_2v_2)$,
$(u_1r_1v_1,r_2)$, $(r_1v_1,u_2r_2)$ and $(u_1r_1,r_2v_2)$ respectively
(see diagram).
In one pass the completion procedure finds all the critical pairs
resulting from overlaps of rules of $R$. Both sides of each of the
critical pairs are reduced as far as possible with respect to $R$ to
obtain a {reduced critical pair} $(c_1,c_2)$. The original
pair
is said to {resolve} if $c_1=c_2$. The reduced pairs that
have
not resolved are orientated, so that $c_1>c_2$, and added to $R$ forming
$R_1$. The procedure is then repeated for the rewrite system $R_1$, to
obtain $R_2$ and so on. When all the critical pairs of a system $R_n$
resolve (i.e. $R_{n+1}=R_n$) then $R_n$ is a {complete rewrite
system}.\\

In terms of Gr\"obner basis theory applied to this special case we
consider the {basis}  $F$
which consists of a {set of two-term polynomials} of the
form
$l-r$ multiplied by $\pm 1$ so that $l>r$.
A {monomial} $m \in X^\dagger$ may be {reduced} 
with respect to $F$ if it contains the {leading monomial} $l$
of a polynomial $l-r$ as a {submonomial} i.e. if $m=ulv$ for
some $u,v \in X^*$.
To reduce $m=ulv$ using the polynomial $l-r$ we replace $l$ by the
{remainder} $r$ of the polynomial, and write $ulv \to_F
urv$.
The Buchberger algorithm looks for {matches between polynomials}.
Given a pair of polynomials $l_1-r_1$, $l_2-r_2$ there are four
possible ways in which an match can occur: $l_1=u_2l_2v_2$,
$u_1l_1v_1=l_2$, $l_1v_1=u_2l_2$ and $u_1l_1=l_2v_2$. The {S-
polynomial} resulting from a match is the difference between the
pair of monomials resulting from applying each two-term polynomial to
the smallest monomial on which the match occurs. The S-polynomials
resulting from each of the four matches are: $r_1-u_2r_2v_2$,
$u_1r_1-v_1,r_2$, $r_1v_1-u_2r_2$ and $u_1r_1-r_2v_2$ respectively
(see diagram). 
In one pass the completion procedure finds all the S-polynomials
resulting from matches of polynomials of $F$. The S-polynomials are
reduced as far as possible with respect to $F$ to obtain a
{reduced
S-polynomial} $c_1-c_2$. Note that reduction can only replace one
term with another so the reduced S-ploynomial will have two terms unless
the two terms reduce to the same thing $c_1=c_2$ in which case the
original S-polynomial is said to {reduce to zero}. The
reduced
S-polynomials that have not been reduced to zero are multiplied by $\pm
1$, so that $c_1>c_2$, and added to $F$ forming $F_1$. The procedure is
then repeated for the basis $F_1$, to obtain $F_2$ and so on. When all
the S-polynomials of a basis $F_n$ reduce to zero (i.e. $F_{n+1}=F_n$)
then $F_n$ is a {Gr\"obner basis}.\\

A critical pair in $R$ will occur if and only if a corresponding 
S-polynomial occurs in $F$. Reduction of the pair by $R$ is equivalent
to reduction of the S-polynomial by $F$. Therefore at any stage any new
rules correspond to the new two-term polynomials and $F_i:=\{l-r:(l,r)
\in R_i\}$.
Therefore the completion procedures as applied to $R$ and $F$ correspond
to each other at every step. {\hfill $\Box$}

\section{Illustration}
This is a picture of the correspondence (viii) and (ix) between critical pairs and 
S-polynomials and the four ways in which they can occur, as described in
the above proof.

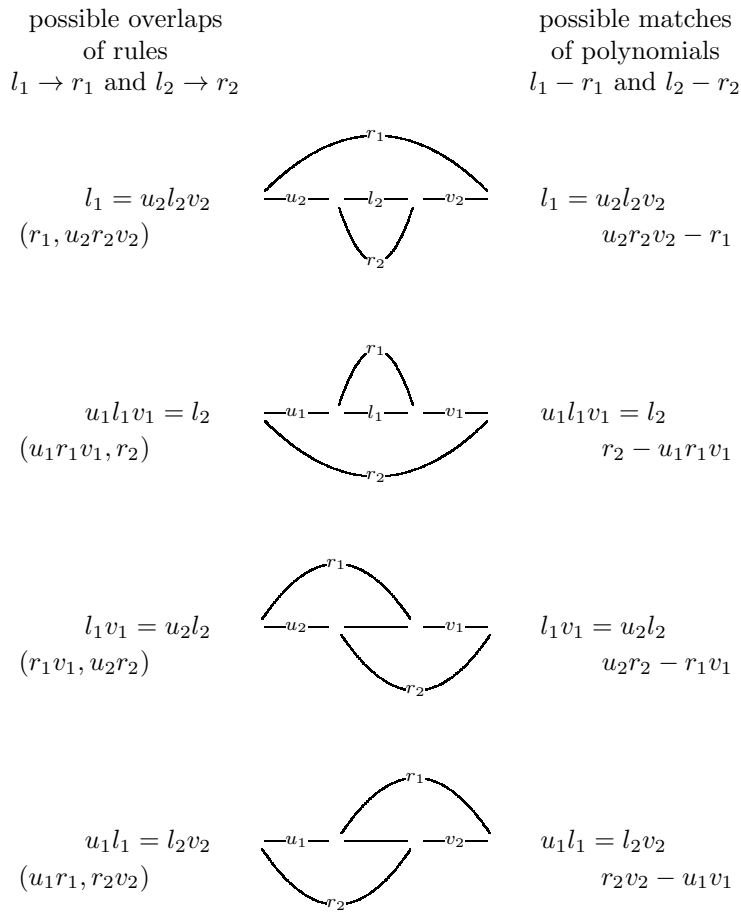
\begin{center}
\begin{tabular}{ccc}
possible overlaps & \mbox{\hspace{3cm}} & possible matches\\
 of rules         &                     &  of polynomials\\
$l_1 \to r_1$ and $l_2 \to r_2$ &   & $l_1 - r_1$ and $l_2 - r_2$\\
\end{tabular}
\end{center}

$$
l_1=u_2l_2v_2 \hspace{0.5cm}
\xymatrix{\ar@{-}@/^2pc/[rrr]|{r_1} \ar@{-}[r]|{u_2} &
\ar@{-}@/_2pc/[r]|{r_2}
\ar@{-}[r]|{l_2} & \ar@{-}[r]|{v_2} & \\}
\hspace{0.5cm}  l_1=u_2l_2v_2
$$

\vspace{-1cm}

$$(r_1,u_2r_2v_2) \makebox[6cm]{}  u_2r_2v_2-r_1$$

\vspace{0.7cm}

$$
u_1l_1v_1=l_2 \hspace{0.5cm}
\xymatrix{\ar@{-}@/_2pc/[rrr]|{r_2} \ar@{-}[r]|{u_1} &
\ar@{-}@/^2pc/[r]|{r_1}
          \ar@{-}[r]|{l_1} & \ar@{-}[r]|{v_1} & \\}
\hspace{0.5cm} u_1l_1v_1=l_2
$$

\vspace{-1cm}

$$(u_1r_1v_1,r_2) \makebox[6cm]{}  r_2-u_1r_1v_1$$

\vspace{0.7cm}

$$
l_1v_1=u_2l_2 \hspace{0.5cm}
\xymatrix{\ar@{-}@/^2pc/[rr]|{r_1} \ar@{-}[r]|{u_2} &
\ar@{-}@/_2pc/[rr]|{r_2}
\ar@{-}[r] & \ar@{-}[r]|{v_1} &\\}
\hspace{0.5cm} l_1v_1=u_2l_2
$$

\vspace{-1cm}

$$
(r_1v_1,u_2r_2)  \makebox[6cm]{} u_2r_2-r_1v_1
$$

\vspace{0.7cm}

$$
u_1l_1=l_2v_2 \hspace{0.5cm}
\xymatrix{\ar@{-}@/_2pc/[rr]|{r_2} \ar@{-}[r]|{u_1} &
\ar@{-}@/^2pc/[rr]|{r_1}
\ar@{-}[r] & \ar@{-}[r]|{v_2} &\\}
\hspace{0.5cm} u_1l_1=l_2v_2
$$

\vspace{-1cm}

$$
(u_1r_1,r_2v_2)  \makebox[6cm]{} r_2v_2-u_1v_1
$$

\vspace {.5cm}

\section{Remarks}
The result that the Knuth-Bendix algorithm is a special case of the 
noncommutative Buchberger algorithm is something that requires further 
investigation. Rewriting techniques and the Knuth-Bendix algorithm have
recently been applied to presentations of Kan extensions over sets \cite{Anne}
and it is not immediately obvious what this will imply for noncommutative 
Gr\"obner bases. Another interesting line of investigation would be to
attempt to adapt rewriting procedures for constructing crossed resolutions of 
group presentations \cite{Anne} to the more general Gr\"obner basis situation.

\textbf{\Large Acknowledgements}\\
This research work was supported by an Earmarked EPSRC Research Studentship
1995-98 `Identities among Relations for Monoids and Categories'.
I would like to thank Larry Lambe for pointing out this likely connection,
Ronnie Brown for discussion on the paper, and Neil Dennis for encouraging my
work.\\

\end{document}